\documentclass[12pt]{article}%
\usepackage{amsfonts}
\usepackage{sw20bams}
\usepackage{amsmath}
\usepackage{amssymb}
\usepackage{graphicx}%
\providecommand{\U}[1]{\protect\rule{.1in}{.1in}}
\begin{document}

\title{Permute, Graph, Map, Derange}
\author{Steven Finch}
\date{January 23, 2022}
\maketitle

\begin{abstract}
We study decomposable combinatorial labeled structures in the exp-log class,
specifically, two examples of type $a=1$ and two examples of type $a=1/2$.
\ Our approach is to establish how well existing theory matches experimental
data. \ For instance, the median length of the longest cycle in a random
$n$-permutation is $(0.6065...)n$, whereas the median length of the largest
component in a random $n$-mapping is $(0.7864...)n$. \ Unsolved problems are
highlighted, in the hope that someone else might address these someday.

\end{abstract}

\footnotetext{Copyright \copyright \ 2021--2022 by Steven R. Finch. All rights
reserved.}Permutations and derangements decompose into cycles; undirected
$2$-regular labeled graphs and mappings decompose into connected components.
\ Among the most striking features of a combinatorial object are

\begin{itemize}
\item the number of cycles or components,

\item the size of the longest cycle or largest component,

\item the size of the shortest cycle or smallest component.
\end{itemize}

\noindent We shall focus on the latter two topics. \ Throughout this paper, a
random object is chosen uniformly from a set, e.g., all
permutations/derangements on $n$ symbols, all $2$-regular graphs with $n$
vertices, and all mappings $\{1,2,\ldots,n\}\rightarrow\{1,2,\ldots,n\}$,
weighted with equal probability. \ Define $b_{n}$ to be the number of
$n$-objects:%
\[
b_{n}=\left\{
\begin{array}
[c]{lll}%
n! &  & \text{for permutations,}\\
\left\vert
\begin{array}
[c]{lll}%
(n-1)b_{n-1}+\tbinom{n-1}{2}b_{n-3} &  & \text{if }n\geq3,\\
b_{0}=1,\;b_{1}=b_{2}=0 &  &
\end{array}
\right.  &  & \text{for graphs,}\\
n^{n} &  & \text{for mappings,}\\
\left\vert
\begin{array}
[c]{lll}%
n\,b_{n-1}+(-1)^{n} &  & \text{if }n\geq1,\\
b_{0}=1 &  &
\end{array}
\right.  &  & \text{for derangements}%
\end{array}
\right.
\]
and $c_{n}$ to be the number of $n$-objects that are connected, i.e., who
possess exactly one component (having $n$ nodes):%
\[
c_{n}=\left\{
\begin{array}
[c]{lll}%
(n-1)! &  & \text{for permutations,}\\
\left\vert
\begin{array}
[c]{lll}%
(n-1)!/2 &  & \text{if }n\geq3,\\
0 &  & \text{otherwise}%
\end{array}
\right.  &  & \text{for graphs,}\\
n!%
{\displaystyle\sum\limits_{j=1}^{n}}
\dfrac{n^{n-j-1}}{(n-j)!} &  & \text{for mappings,}\\
\left\vert
\begin{array}
[c]{lll}%
(n-1)! &  & \text{if }n\geq2,\\
0 &  & \text{otherwise}%
\end{array}
\right.  &  & \text{for derangements.}%
\end{array}
\right.
\]

Our starting point is a recursive formula for $L_{k,n}$, the number of
$n$-objects whose largest component has exactly $k$ nodes, $1\leq k\leq n$.
\ The initial conditions are%
\[%
\begin{array}
[c]{ccc}%
L_{0,n}=\delta_{0,n}, &  & L_{1,n}=1-\delta_{0,n}%
\end{array}
\]
for permutations and mappings, and%
\[%
\begin{array}
[c]{ccc}%
L_{0,n}=\delta_{0,n}, &  & L_{1,n}=0
\end{array}
\]
for graphs and derangements. \ Letting%
\[
m_{j,k,n}=\min\{k-1,n-k\,j\}
\]
(although we suppress the dependence on $j$, $k$, $n$), we have%
\[
L_{k,n}=%
{\displaystyle\sum\limits_{j=1}^{\left\lfloor n/k\right\rfloor }}
\frac{n!c_{k}^{j}}{j!(k!)^{j}(n-k\,j)!}\cdot\left\{
\begin{array}
[c]{ccc}%
{\displaystyle\sum\limits_{i=1}^{m}}
L_{i,n-k\,j} &  & \text{if }m\geq1,\\%
{\displaystyle\sum\limits_{i=0}^{1}}
L_{i,n-k\,j} &  & \text{if }m=0.
\end{array}
\right.
\]
The separation of cases $m\geq1$ and $m=0$ clarifies what was surely intended
in \cite{GG-tcs6, PR1-tcs6}. \ 

Next is a recursive formula for $S_{k,n}$, the number of $n$-objects whose
smallest component has exactly $k$ nodes, $1\leq k\leq n$. \ The initial
conditions are%
\[%
\begin{array}
[c]{c}%
S_{0,n}=\delta_{0,n}%
\end{array}
\]
for permutations and mappings, and%
\[%
\begin{array}
[c]{c}%
S_{0,n}=0
\end{array}
\]
for graphs and derangements. \ Letting%
\[
\theta_{k,n}=\left\{
\begin{array}
[c]{lll}%
1 &  & \text{if }k\text{ is a divisor of }n\text{,}\\
0 &  & \text{otherwise}%
\end{array}
\right.
\]
we have \cite{PR1-tcs6}%
\[
S_{k,n}=%
{\displaystyle\sum\limits_{j=1}^{\left\lfloor n/k\right\rfloor }}
\frac{n!c_{k}^{j}}{j!(k!)^{j}(n-k\,j)!}\cdot%
{\displaystyle\sum\limits_{i=k+1}^{n-k\,j}}
S_{i,n-k\,j}+\theta_{k,n}\frac{n!c_{k}^{n/k}}{(n/k)!(k!)^{n/k}}.
\]
Clearly
\[%
{\displaystyle\sum\limits_{k=1}^{n}}
L_{k,n}=%
{\displaystyle\sum\limits_{k=1}^{n}}
S_{k,n}=b_{n}%
\]
and $L_{n,n}=S_{n,n}=c_{n}$ (the latter following from the empty sum
convention). \ A computer algebra software package (e.g., Mathematica) makes
exact integer calculations for ample $n$ of $L_{k,n}$ and $S_{k,n}$ feasible.

Permutations and derangements belong to the exp-log class of type $a=1$,
whereas graphs and mappings belong to the exp-log class of type $a=1/2$.
\ Explaining the significance of the parameter $a>0$ would take us too far
afield \cite{PR2-tcs6}. \ Let%
\[%
\begin{array}
[c]{ccc}%
E(x)=%
{\displaystyle\int\limits_{x}^{\infty}}
\dfrac{e^{-t}}{t}dt=-\operatorname{Ei}(-x), &  & x>0
\end{array}
\]
be the exponential integral. \ Define \cite{SL-tcs6, Shi-tcs6, FO-tcs6,
Gou-tcs6, ABT-tcs6, Pin-tcs6}%
\[
_{L}G_{a}(r,h)=\frac{\Gamma(a+1)a^{r-1}}{\Gamma(a+h)(r-1)!}%
{\displaystyle\int\limits_{0}^{\infty}}
x^{h-1}E(x)^{r-1}\exp\left[  -a\,E(x)-x\right]  dx,
\]%
\[
_{S}G_{a}(r,h)=\left\{
\begin{array}
[c]{lll}%
e^{-h\,\gamma}a^{r-1}/r! &  & \text{if }h=a,\\
\dfrac{\Gamma(a+1)}{(h-1)!(r-1)!}%
{\displaystyle\int\limits_{0}^{\infty}}
x^{h-1}\exp\left[  a\,E(x)-x\right]  dx &  & \text{if }h>a
\end{array}
\right.
\]
which are related to the $h^{\text{th}}$ moment of the $r^{\text{th}}$
largest/smallest component size (in this paper, rank $r=1$; height $h=1$ or
$2$). \ Our notation $_{S}G_{a}$ is deceiving. \ While permutation
and\ derangement moments coincide for $L$ (both being $_{L}G_{a}$ with $a=1$),
they are \textbf{not} equal for $S$ (they differ by a factor $e$). \ In the
same way, graph and\ mapping moments coincide for $L$ but differ for $S$ (by a
factor $e^{3/4}/\sqrt{2}$).

Finally, we ask the first (of several)\ questions. \ Does a general recursion
(similar to that for $L_{k,n}$ and $S_{k,n}$) exist for the number $N_{k,n}$
of $n$ objects who possess exactly $k$ components? \ Specific formulas are
known for each of our four examples, but they are quite dissimilar. \ It seems
unlikely that generalization to arbitrary $N_{k,n}$ is possible, however it is
still well-worth contemplating.

\section{Permute}

For fixed $n$, the sequences $\{L_{k,n}:1\leq k\leq n\}$ and $\{S_{k,n}:1\leq
k\leq n\}$ constitute probability mass functions (upon normalization by
$b_{n}=n!$) \ These have corresponding means $_{L}\mu_{n}$, $_{S}\mu_{n}$ and
variances $_{L}\sigma_{n}^{2}$, $_{S}\sigma_{n}^{2}$ given in Table 1. \ We
also provide the median $_{L}\nu_{n}$; note that $_{S}\nu_{n}=1$ for $n\geq3$
is trivial. \ For convenience (in table headings only), the following notation
is used:
\[%
\begin{array}
[c]{ccccc}%
_{L}\widetilde{\mu}_{n}=\dfrac{_{L}\mu_{n}}{n}, &  & _{L}\widetilde{\sigma
}_{n}^{2}=\dfrac{_{L}\sigma_{n}^{2}}{n^{2}}, &  & _{L}\widetilde{\nu}%
_{n}=\dfrac{_{L}\nu_{n}}{n},
\end{array}
\]%
\[%
\begin{array}
[c]{ccc}%
_{S}\widetilde{\mu}_{n}=\left\{
\begin{array}
[c]{lll}%
\dfrac{_{S}\mu_{n}}{\ln(n)} & \bigskip & \text{if }a=1,\\
\dfrac{_{S}\mu_{n}}{n^{1/2}} &  & \text{if }a=1/2;
\end{array}
\right.  &  & _{S}\widetilde{\sigma}_{n}^{2}=\left\{
\begin{array}
[c]{lll}%
\dfrac{_{S}\sigma_{n}^{2}}{n} & \bigskip & \text{if }a=1,\\
\dfrac{_{S}\sigma_{n}^{2}}{n^{3/2}} &  & \text{if }a=1/2.
\end{array}
\right.
\end{array}
\]
\medskip

\begin{center}%
\begin{tabular}
[c]{|c|c|c|c|c|c|}\hline
$n$ & $_{L}\widetilde{\mu}_{n}$ & $_{L}\widetilde{\sigma}_{n}^{2}$ &
$_{L}\widetilde{\nu}_{n}$ & $_{S}\widetilde{\mu}_{n}$ & $_{S}\widetilde
{\sigma}_{n}^{2}$\\\hline
1000 & 0.624642 & 0.036945 & 0.6060 & 0.717352 & 1.307043\\\hline
2000 & 0.624486 & 0.036926 & 0.6060 & 0.703135 & 1.307125\\\hline
3000 & 0.624434 & 0.036920 & 0.6063 & 0.695960 & 1.307153\\\hline
4000 & 0.624408 & 0.036917 & 0.6062 & 0.691295 & 1.307167\\\hline
\end{tabular}

Table 1:\ Statistics for Permute ($a=1$)
\end{center}

\noindent We have%
\[
\lim_{n\rightarrow\infty}\dfrac{_{L}\mu_{n}}{n}=\,_{L}G_{1}%
(1,1)=0.62432998854355087099...,
\]%
\[
\lim_{n\rightarrow\infty}\dfrac{_{L}\sigma_{n}^{2}}{n^{2}}=\,_{L}%
G_{1}(1,2)-\,_{L}G_{1}(1,1)^{2}=0.03690783006485220217...,
\]%
\[
\lim_{n\rightarrow\infty}\dfrac{_{L}\nu_{n}}{n}=\frac{1}{\sqrt{e}%
}=0.60653065971263342360...,
\]%
\[
\lim_{n\rightarrow\infty}\dfrac{_{S}\mu_{n}}{\ln(n)}=e^{-\gamma}%
=0.56145948356688516982...,
\]%
\[
\lim_{n\rightarrow\infty}\dfrac{_{S}\sigma_{n}^{2}}{n}=\,_{S}G_{1}%
(1,2)=1.30720779891056809974....
\]

\noindent It is not surprising that $_{S}\sigma_{n}^{2}$\ enjoys linear
growth:\ $S_{1,n}\sim(1-1/e)n!$ and $S_{n,n}=(n-1)!$ jointly place
considerable weight on the distributional extremes. \ The unusual logarithmic
growth of $_{S}\mu_{n}$\ is due to $S_{1,n}$ nevertheless overwhelming all
other $S_{k,n}$.

A one-line proof of the $_{L}\nu_{n}$\ result is \cite{Bz-tcs6, Wk-tcs6}%
\[
\lim_{n\rightarrow\infty}%
{\displaystyle\sum\limits_{k=\left\lfloor \frac{n}{\sqrt{e}}\right\rfloor
}^{n}}
\frac{1}{k}=\lim_{n\rightarrow\infty}\ln(n)-\ln\left(  \frac{n}{\sqrt{e}%
}\right)  =\ln\left(  \frac{1}{\sqrt{e}}\right)  =\frac{1}{2}.
\]
Alternatively, the asymptotic probability that the longest cycle has size
$>n\,x$ is \cite{GG-tcs6}%
\[%
{\displaystyle\int\limits_{x}^{1}}
\frac{1}{y}dy=\frac{1}{2}%
\]
hence $\ln(1/x)=1/2$ and $x=1/\sqrt{e}$. \ We will see a variation of this
approach later.

No formula for the covariance between sizes of the longest cycle and shortest
cycle is known. \ Interplay between the number of cycles and either of the
extremes likewise remains inscrutable. \ We earlier examined not permutations,
but instead integer compositions, finding a complicated recursion for a
certain bivariate probability distribution \cite{NS-tcs6, Fi-tcs6}. \ Thus a
cross-correlation can be estimated for compositions, but not yet for
permutations \cite{O1-tcs6}.

\section{Graph}

Let us explain why $L_{3,6}=10$, $L_{6,6}=60$; $L_{4,7}=105$, $L_{7,7}=360$;
and $L_{4,8}=315$, $L_{5,8}=672$, $L_{8,8}=2520$. \ When $n=6$, a $2$-regular
graph is either a hexagon, with $5!/2$ distinct labelings, or the disjoint
union of two triangles, with $\tbinom{6}{3}/2$ labelings. \ When $n=7$, a
$2$-regular graph is either a heptagon, with $6!/2$ distinct labelings, or the
disjoint union of a triangle and a square, with $3\tbinom{7}{3}$ labelings
(since $3!/2=3$). \ When $n=8$, a $2$-regular graph is either an octagon, with
$7!/2$ distinct labelings; the disjoint union of a triangle and a pentagon,
with $12\tbinom{7}{3}$ labelings (since $4!/2=12$); or the disjoint union of
two squares, with $9\tbinom{8}{4}/2$ labelings (since $(3!/2)^{2}=9$).
\ Circumstances become more complicated when $n=9$: a $2$-regular graph is
either an enneagon, or the disjoint union of a square and a pentagon, or the
disjoint union of three triangles, or the disjoint union of a triangle and a
hexagon \cite{O2-tcs6}.

Upon normalization by $b_{n}$, we obtain

\begin{center}%
\begin{tabular}
[c]{|c|c|c|c|c|c|}\hline
$n$ & $_{L}\widetilde{\mu}_{n}$ & $_{L}\widetilde{\sigma}_{n}^{2}$ &
$_{L}\widetilde{\nu}_{n}$ & $_{S}\widetilde{\mu}_{n}$ & $_{S}\widetilde
{\sigma}_{n}^{2}$\\\hline
1000 & 0.758771 & 0.037099 & 0.7860 & 3.007677 & 2.097084\\\hline
2000 & 0.758297 & 0.037053 & 0.7865 & 3.029960 & 2.096470\\\hline
3000 & 0.758139 & 0.037038 & 0.7863 & 3.039930 & 2.096262\\\hline
4000 & 0.758060 & 0.037030 & 0.7865 & 3.045902 & 2.096157\\\hline
\end{tabular}

Table 2:\ Statistics for Graph ($a=1/2$)
\end{center}

\noindent and%

\[
\lim_{n\rightarrow\infty}\dfrac{_{L}\mu_{n}}{n}=\,_{L}G_{1/2}%
(1,1)=0.75782301126849283774...,
\]%
\[
\lim_{n\rightarrow\infty}\dfrac{_{L}\sigma_{n}^{2}}{n^{2}}=\,_{L}%
G_{1/2}(1,2)-\,_{L}G_{1/2}(1,1)^{2}=0.03700721658229030320...,
\]%
\[
\lim_{n\rightarrow\infty}\dfrac{_{L}\nu_{n}}{n}=\frac{4e}{(1+e)^{2}%
}=0.78644773296592741014...,
\]%
\[
\lim_{n\rightarrow\infty}\dfrac{_{S}\mu_{n}}{n^{1/2}}=e^{3/4}\,_{S}%
G_{1/2}(1,1)=3.08504247563149222958...,
\]%
\[
\lim_{n\rightarrow\infty}\dfrac{_{S}\sigma_{n}^{2}}{n^{3/2}}=e^{3/4}%
\,_{S}G_{1/2}(1,2)=2.09583743942571712967....
\]

\noindent Some non-explicit formulas for the latter two results arise in
Section 7. \ A proof for the median result is deferred to Section 3.

\section{Map}

Let us explain why $S_{1,2}=1$, $S_{2,2}=3$ and $S_{1,3}=10$, $S_{3,3}=17$.
\ The unique $2$-map with totally disconnected nodes is the identity map; we
associate this map with its image sequence $12$, i.e., two isolated loops (two
components of size $1$). \ The maps $11$ and $22$ are each pictured as one
loop attached to a $1$-tail (a component of size $2$); the map $21$ is
pictured as a $2$-cycle (again, a component of size $2$). \ For $3$-maps, we
have%
\[%
\begin{array}
[c]{ccc}%
123 &  & \text{i.e., three isolated loops;}%
\end{array}
\]%
\[%
\begin{array}
[c]{ccc}%
113,121,122,133,223,323 &  & \text{i.e., one isolated loop and one loop
attached to a }1\text{-tail;}%
\end{array}
\]%
\[%
\begin{array}
[c]{ccc}%
213,321,132 &  & \text{i.e., one isolated loop and one }2\text{-cycle}%
\end{array}
\]
which give ten cases, and%
\[%
\begin{array}
[c]{ccc}%
112,131,221,322,233,313 &  & \text{i.e., one loop attached to a }%
2\text{-tail;}%
\end{array}
\]%
\[%
\begin{array}
[c]{ccc}%
111,222,333 &  & \text{i.e., one loop attached to two }1\text{-tails;}%
\end{array}
\]%
\[%
\begin{array}
[c]{ccc}%
211,212,232,311,331,332 &  & \text{i.e., one }2\text{-cycle attached to a
}1\text{-tail;}%
\end{array}
\]%
\[%
\begin{array}
[c]{ccc}%
231,312 &  & \text{i.e., one }3\text{-cycle}%
\end{array}
\]
which give seventeen cases \cite{SF-tcs6, O3-tcs6}. \ The complexity grows
when $n=4$: it can be shown that $S_{1,4}=87$, $S_{2,4}=27$, $S_{4,4}=142$.

Upon normalization by $b_{n}=n^{n}$, we obtain

\begin{center}%
\begin{tabular}
[c]{|c|c|c|c|c|c|}\hline
$n$ & $_{L}\widetilde{\mu}_{n}$ & $_{L}\widetilde{\sigma}_{n}^{2}$ &
$_{L}\widetilde{\nu}_{n}$ & $_{S}\widetilde{\mu}_{n}$ & $_{S}\widetilde
{\sigma}_{n}^{2}$\\\hline
1000 & 0.762505 & 0.036968 & 0.7920 & 1.969526 & 1.384968\\\hline
2000 & 0.761122 & 0.036980 & 0.7905 & 1.991932 & 1.389355\\\hline
3000 & 0.760512 & 0.036985 & 0.7900 & 2.002505 & 1.391309\\\hline
4000 & 0.760149 & 0.036988 & 0.7895 & 2.009048 & 1.392477\\\hline
\end{tabular}

Table 3:\ Statistics for Map ($a=1/2$)
\end{center}

\noindent and%

\[
\lim_{n\rightarrow\infty}\dfrac{_{L}\mu_{n}}{n}=\,_{L}G_{1/2}%
(1,1)=0.75782301126849283774...,
\]%
\[
\lim_{n\rightarrow\infty}\dfrac{_{L}\sigma_{n}^{2}}{n^{2}}=\,_{L}%
G_{1/2}(1,2)-\,_{L}G_{1/2}(1,1)^{2}=0.03700721658229030320...,
\]%
\[
\lim_{n\rightarrow\infty}\dfrac{_{L}\nu_{n}}{n}=\frac{4e}{(1+e)^{2}%
}=0.78644773296592741014...,
\]%
\[
\lim_{n\rightarrow\infty}\dfrac{_{S}\mu_{n}}{n^{1/2}}=\sqrt{2}\,_{S}%
G_{1/2}(1,1)=2.06089224152016653900...,
\]%
\[
\lim_{n\rightarrow\infty}\dfrac{_{S}\sigma_{n}^{2}}{n^{3/2}}=\sqrt{2}%
\,_{S}G_{1/2}(1,2)=1.40007638550124502818....
\]

\noindent As before, some non-explicit formulas for the latter two results
arise in Section 7. \ The asymptotic probability that the largest component
has size $>n\,x$ is \cite{Kol-tcs6, DEP-tcs6}%
\[%
{\displaystyle\int\limits_{x}^{1}}
\frac{1}{2y\sqrt{1-y}}dy=\frac{1}{2}%
\]
hence%
\[
\frac{1}{2}\ln\left(  \frac{1+\sqrt{1-x}}{1-\sqrt{1-x}}\right)  =\frac{1}{2}%
\]
and $x=4e/(1+e)^{2}$. \ 

\section{Derange}

Derangements are permutations with no fixed points \cite{O4-tcs6}. \ It is
easy to show that $L_{3,5}=20$ (a longest cycle in a $5$-derangement cannot
have size $2$ or $4$) and $S_{2,5}=20$. \ Upon normalization by $b_{n}$, we obtain

\begin{center}%
\begin{tabular}
[c]{|c|c|c|c|c|c|}\hline
$n$ & $_{L}\widetilde{\mu}_{n}$ & $_{L}\widetilde{\sigma}_{n}^{2}$ &
$_{L}\widetilde{\nu}_{n}$ & $_{S}\widetilde{\mu}_{n}$ & $_{S}\widetilde
{\sigma}_{n}^{2}$\\\hline
1000 & 0.625266 & 0.037018 & 0.6060 & 1.701217 & 3.551193\\\hline
2000 & 0.624798 & 0.036963 & 0.6065 & 1.685257 & 3.552276\\\hline
3000 & 0.624642 & 0.036945 & 0.6067 & 1.677202 & 3.552637\\\hline
4000 & 0.624564 & 0.036935 & 0.6065 & 1.671965 & 3.552818\\\hline
\end{tabular}

Table 4:\ Statistics for Derange ($a=1$)
\end{center}

\noindent and%
\[
\lim_{n\rightarrow\infty}\dfrac{_{L}\mu_{n}}{n}=\,_{L}G_{1}%
(1,1)=0.62432998854355087099...,
\]%
\[
\lim_{n\rightarrow\infty}\dfrac{_{L}\sigma_{n}^{2}}{n^{2}}=\,_{L}%
G_{1}(1,2)-\,_{L}G_{1}(1,1)^{2}=0.03690783006485220217...,
\]%
\[
\lim_{n\rightarrow\infty}\dfrac{_{L}\nu_{n}}{n}=\frac{1}{\sqrt{e}%
}=0.60653065971263342360...,
\]

\[
\lim_{n\rightarrow\infty}\dfrac{_{S}\mu_{n}}{\ln(n)}=e^{1-\gamma
}=1.52620511159586388047...,
\]%
\[
\lim_{n\rightarrow\infty}\dfrac{_{S}\sigma_{n}^{2}}{n}=e\cdot\,_{S}%
G_{1}(1,2)=3.55335920579854297440....
\]

\noindent The asymptotic expression for the average shortest cycle length
follows from
\[
\frac{_{S}\mu_{n}\cdot b_{n}+1\cdot\left(  n!-b_{n}\right)  }{n!}\sim
e^{-\gamma}\ln(n)
\]
and the fact that $b_{n}/n!\rightarrow1/e$ as $n\rightarrow\infty$; similarly
for higher moments.

\section{Generalized Dickman Rho (I)}

Define $b_{n,m}$ to be the number of $n$-objects whose largest component has
size $\leq m$; thus $b_{n}=b_{n,n}$. \ Given $a>0$, let%
\[
\rho_{a}(x)=\left\{
\begin{array}
[c]{lll}%
1 &  & \text{if }0\leq x<1,\\
1-a%
{\displaystyle\int\limits_{1}^{x}}
\dfrac{\rho_{a}(t-1)\cdot(t-1)^{a-1}}{t^{a}}dt &  & \text{if }x\geq1
\end{array}
\right.
\]
and observe that the standard Dickman function $\rho(x)=\rho_{1}(x)$. \ A
theorem proved in \cite{OPRW-tcs6} asserts that%
\[
\lim_{m\rightarrow\infty}\frac{b_{\left\lfloor x\,m\right\rfloor ,m}%
}{b_{\left\lfloor x\,m\right\rfloor }}=\rho_{a}(x)
\]
for any $x>1$. \ Of course,%
\[%
{\displaystyle\sum\limits_{j=1}^{m}}
L_{j,n}=b_{n,m}%
\]
hence we can easily verify this result experimentally. \ 

\begin{center}%
\begin{tabular}
[c]{|l|l|l|l|l|l|l|l|l|}\hline
$m\setminus x$ & 2 & 3 & 4 & 5 &  & 2 & 3 & 4\\\hline
100 & 0.309347 & 0.049634 & 0.0050952 & 0.0003748 &  & 0.117715 & 0.0082644 &
0.0003680\\\hline
200 & 0.308101 & 0.049121 & 0.0050026 & 0.0003646 &  & 0.118178 & 0.0082399 &
0.0003638\\\hline
300 & 0.307685 & 0.048950 & 0.0049719 & 0.0003613 &  & 0.118329 & 0.0082309 &
0.0003624\\\hline
400 & 0.307477 & 0.048864 & 0.0049566 & 0.0003597 &  & 0.118404 & 0.0082262 &
0.0003616\\\hline
500 & 0.307353 & 0.048813 & 0.0049475 & 0.0003587 &  & 0.118449 & 0.0082233 &
0.0003612\\\hline
600 & 0.307269 & 0.048779 & 0.0049414 & 0.0003580 &  & 0.118478 & 0.0082214 &
0.0003609\\\hline
700 & 0.307210 & 0.048755 & 0.0049370 & 0.0003575 &  & 0.118500 & 0.0082200 &
0.0003607\\\hline
800 & 0.307165 & 0.048736 & 0.0049337 & 0.0003572 &  & 0.118516 & 0.0082190 &
0.0003605\\\hline
$\vdots$ & $\vdots$ & $\vdots$ & $\vdots$ & $\vdots$ &  & $\vdots$ & $\vdots$
& $\vdots$\\\hline
$\infty$ & 0.306853 & 0.048608 & 0.0049109 & 0.0003547 &  & 0.118626 &
0.0082115 & 0.0003594\\\hline
\end{tabular}

Table 5A:\ Ratio $b_{\left\lfloor x\,m\right\rfloor ,m}/b_{\left\lfloor
x\,m\right\rfloor }$ for Permute ($a=1$) and Graph ($a=1/2$)
\end{center}

\smallskip

\begin{center}%
\begin{tabular}
[c]{|l|l|l|l|l|l|l|l|l|}\hline
$m\setminus x$ & 2 & 3 & 4 &  & 2 & 3 & 4 & 5\\\hline
100 & 0.111305 & 0.0074576 & 0.0003185 &  & 0.304359 & 0.048597 & 0.0049699 &
0.0003645\\\hline
200 & 0.112756 & 0.0076060 & 0.0003258 &  & 0.305604 & 0.048605 & 0.0049409 &
0.0003596\\\hline
300 & 0.113579 & 0.0076901 & 0.0003302 &  & 0.306020 & 0.048607 & 0.0049310 &
0.0003580\\\hline
400 & 0.114124 & 0.0077458 & 0.0003332 &  & 0.306228 & 0.048608 & 0.0049260 &
0.0003572\\\hline
500 & 0.114518 & 0.0077862 & 0.0003354 &  & 0.306353 & 0.048608 & 0.0049230 &
0.0003567\\\hline
600 & 0.114822 & 0.0078173 & 0.0003371 &  & 0.306436 & 0.048608 & 0.0049210 &
0.0003564\\\hline
700 & 0.115065 & 0.0078423 & 0.0003384 &  & 0.306496 & 0.048608 & 0.0049196 &
0.0003561\\\hline
800 & 0.115265 & 0.0078629 & 0.0003396 &  & 0.306540 & 0.048608 & 0.0049185 &
0.0003559\\\hline
$\vdots$ & $\vdots$ & $\vdots$ & $\vdots$ &  & $\vdots$ & $\vdots$ & $\vdots$
& $\vdots$\\\hline
$\infty$ & 0.118626 & 0.0082115 & 0.0003594 &  & 0.306853 & 0.048608 &
0.0049109 & 0.0003547\\\hline
\end{tabular}

Table 5B:\ Ratio $b_{\left\lfloor x\,m\right\rfloor ,m}/b_{\left\lfloor
x\,m\right\rfloor }$ for Map ($a=1/2$) and Derange ($a=1$)
\end{center}

Why have we devoted effort to evaluating Dickman's rho? \ Answer:\ the
function $\rho_{a}(x)$ is fundamentally connected to $L_{k,n}$ asymptotics in
Sections 1--4. \ The $h^{\text{th}}$ moments of the largest component size are%
\[%
\begin{array}
[c]{ccccccc}%
{\displaystyle\int\limits_{0}^{\infty}}
\dfrac{\rho_{1}(x)}{(x+1)^{h+1}}dx &  & \text{and} &  & \dfrac{1}{2}%
{\displaystyle\int\limits_{0}^{\infty}}
\dfrac{\rho_{1/2}(x)}{\sqrt{x}(x+1)^{h+1/2}}dx, &  & h=1,2,\ldots
\end{array}
\]
for $a=1$ and $a=1/2$, respectively. \ Extension to arbitrary $a>0$ is
possible. \ Of course, we also have integrals $_{L}G_{a}(r,h)$ available.

\section{Generalized Buchstab Omega (I)}

Define $b_{n,m}$ to be the number of $n$-objects whose smallest component has
size $\geq m$; note that $c_{n}=b_{n,n}$. \ Given $a>0$, let%
\[
\Omega_{a}(x)=\left\{
\begin{array}
[c]{lll}%
1 &  & \text{if }1\leq x<2,\\
1+a%
{\displaystyle\int\limits_{2}^{x}}
\dfrac{\Omega_{a}(t-1)}{t-1}dt &  & \text{if }x\geq2
\end{array}
\right.
\]
and observe that the standard Buchstab function $\omega(x)=\Omega_{1}(x)/x$.
\ A theorem proved in \cite{BMPR-tcs6} asserts that%
\[
\lim_{m\rightarrow\infty}\frac{c_{\left\lfloor x\,m\right\rfloor }%
}{b_{\left\lfloor x\,m\right\rfloor ,m}}=\frac{1}{\Omega_{a}(x)}%
\]
for any $x>1$. \ Of course,%
\[%
{\displaystyle\sum\limits_{j=m}^{n}}
S_{j,n}=b_{n,m}%
\]
hence we can easily verify this result experimentally.

As an aside, $c_{n}/b_{n,m}$ is called the \textit{probability of
connectedness} in \cite{BMPR-tcs6}, i.e., the odds that an $n$-object, whose
smallest component has size at least $m$, is connected. \ No analogous name
has been proposed for $b_{n,m}/b_{n}$ from Section 5, i.e., the odds that all
components of an $n$-object have size at most $m$. \ Maybe \textit{probability
of smoothness} would be appropriate (\textquotedblleft
smooth\textquotedblright\ coming from prime number theory). \ For Section 7,
the same ratio might be called the \textit{probability of roughness}, wherein
all components of an $n$-object have size at least $m$.

\begin{center}%
\begin{tabular}
[c]{|l|l|l|l|l|l|l|l|l|l|}\hline
$m\setminus x$ & 2 & 3 & 4 & 5 &  & 2 & 3 & 4 & 5\\\hline
100 & 0.990 & 0.587992 & 0.443034 & 0.354438 &  & 0.995 & 0.740555 &
0.628689 & 0.555092\\\hline
200 & 0.995 & 0.589306 & 0.444151 & 0.355327 &  & 0.997 & 0.741591 &
0.629581 & 0.555860\\\hline
300 & 0.997 & 0.589743 & 0.444523 & 0.355624 &  & 0.998 & 0.741936 &
0.629878 & 0.556116\\\hline
400 & 0.997 & 0.589962 & 0.444710 & 0.355772 &  & 0.999 & 0.742108 &
0.630027 & 0.556244\\\hline
500 & 0.998 & 0.590093 & 0.444822 & 0.355862 &  & 0.999 & 0.742212 &
0.630116 & 0.556321\\\hline
600 & 0.998 & 0.590180 & 0.444896 & 0.355921 &  & 0.999 & 0.742281 &
0.630175 & 0.556372\\\hline
700 & 0.999 & 0.590242 & 0.444949 & 0.355963 &  & 0.999 & 0.742330 &
0.630218 & 0.556409\\\hline
800 & 0.999 & 0.590289 & 0.444989 & 0.355995 &  & 0.999 & 0.742367 &
0.630250 & 0.556436\\\hline
$\vdots$ & $\vdots$ & $\vdots$ & $\vdots$ & $\vdots$ &  & $\vdots$ & $\vdots$
& $\vdots$ & $\vdots$\\\hline
$\infty$ & 1 & 0.590616 & 0.445269 & 0.356218 &  & 1 & 0.742626 & 0.630473 &
0.556628\\\hline
\end{tabular}

Table 6A:\ Ratio $c_{\left\lfloor x\,m\right\rfloor }/b_{\left\lfloor
x\,m\right\rfloor ,m}$ for Permute ($a=1$) and Graph ($a=1/2$)
\end{center}

\smallskip

\begin{center}%
\begin{tabular}
[c]{|l|l|l|l|l|l|l|l|l|l|}\hline
$m\setminus x$ & 2 & 3 & 4 & 5 &  & 2 & 3 & 4 & 5\\\hline
100 & 0.995 & 0.746112 & 0.635215 & 0.561960 &  & 0.990 & 0.587992 &
0.443034 & 0.354438\\\hline
200 & 0.998 & 0.745502 & 0.634175 & 0.560698 &  & 0.995 & 0.589306 &
0.444151 & 0.355327\\\hline
300 & 0.998 & 0.745124 & 0.633623 & 0.560060 &  & 0.997 & 0.589743 &
0.444523 & 0.355624\\\hline
400 & 0.999 & 0.744866 & 0.633267 & 0.559656 &  & 0.997 & 0.589962 &
0.444710 & 0.355772\\\hline
500 & 0.999 & 0.744677 & 0.633013 & 0.559371 &  & 0.998 & 0.590093 &
0.444822 & 0.355862\\\hline
600 & 0.999 & 0.744530 & 0.632819 & 0.559156 &  & 0.998 & 0.590180 &
0.444896 & 0.355921\\\hline
700 & 0.999 & 0.744412 & 0.632664 & 0.558985 &  & 0.999 & 0.590242 &
0.444949 & 0.355963\\\hline
800 & 0.999 & 0.744314 & 0.632537 & 0.558845 &  & 0.999 & 0.590289 &
0.444989 & 0.355995\\\hline
$\vdots$ & $\vdots$ & $\vdots$ & $\vdots$ & $\vdots$ &  & $\vdots$ & $\vdots$
& $\vdots$ & $\vdots$\\\hline
$\infty$ & 1 & 0.742626 & 0.630473 & 0.556628 &  & 1 & 0.590616 & 0.445269 &
0.356218\\\hline
\end{tabular}

Table 6B:\ Ratio $c_{\left\lfloor x\,m\right\rfloor }/b_{\left\lfloor
x\,m\right\rfloor ,m}$ for Map ($a=1/2$) and Derange ($a=1$)
\end{center}

Buchstab's Omega, as defined here, does not seem to be allied with $S_{k,n}$
asymptotics in Sections 1--4. \ A different generalization is discussed in
Section 7.

\section{Generalized Buchstab Omega (II)}

Define $b_{n,m}$ to be the number of $n$-objects whose smallest component has
size $\geq m$ (as in Section 6). \ When restricting attention to permutations,
Panario \& Richmond \cite{PR2-tcs6} obtained that
\[
\lim_{m\rightarrow\infty}\frac{m\,b_{\left\lfloor x\,m\right\rfloor ,m}%
}{b_{\left\lfloor x\,m\right\rfloor }}=\omega(x)
\]
for any $x>1$, where $\omega(x)$ is the standard Buchstab function. \ They
seemed to presume that the same limit would occur for derangements (since both
permutations and derangements have $a=1$), which is not true. \ Replace now
the initial factor $m$ in the numerator by $m^{1/2}$. \ Panario \& Richmond
realized that $2$-regular graphs and mappings would possess a limit different
from $\omega(x)$. \ They seemed, however, to presume that equivalent limits
would occur (since both graphs and maps have $a=1/2$), which is again untrue.
\ In Section 5, we studied two functions $\Omega_{a}$, $a\in\{1,1/2\}$; here
we have four omega (lowercase \textquotedblleft o\textquotedblright) functions
$\omega_{A}$, $A\in\{P,G,M,D\}$, one for each structure under consideration.
\ Upon multiplication of limits, we discover%
\[
\frac{\omega_{A}(x)}{\Omega_{a}(x)}=\lim_{m\rightarrow\infty}\frac
{m^{a}c_{\left\lfloor x\,m\right\rfloor }}{b_{\left\lfloor x\,m\right\rfloor
}}=\frac{1}{x^{a}}\cdot\left\{
\begin{array}
[c]{lll}%
1 &  & \text{if }A=P\text{ and }a=1,\\
e^{3/4}\sqrt{\pi}/2 &  & \text{if }A=G\text{ and }a=1/2,\\
\sqrt{\pi/2} &  & \text{if }A=M\text{ and }a=1/2,\\
e &  & \text{if }A=D\text{ and }a=1
\end{array}
\right.
\]
using known $c_{n}/b_{n}$ asymptotics as $n\rightarrow\infty$ for graphs and
mappings \cite{Ktz-tcs6, Cmt-tcs6, O5-tcs6}. \ Perhaps, for fixed $a>0$,
$\omega_{A}$ varies only up to multiplicative constant. \ These formulas allow
us to provide numerical values in the final rows of Tables 6A\ and 6B.

Return now to Panario \& Richmond. Especially puzzling is a claim (for
permutations)\ that \cite{PR2-tcs6}
\[
\lim_{n\rightarrow\infty}\dfrac{_{S}\sigma_{n}^{2}}{n}=%
{\displaystyle\int\limits_{2}^{\infty}}
\dfrac{\omega(x)}{x^{2}}dx.
\]
From Section 1, the left-hand side is $_{S}G_{1}(1,2)=1.307207...$
\cite{SL-tcs6} whereas the right-hand side is $\frac{1}{2}(0.556816...)$
\cite{PR3-tcs6}. \ Thus predictions in \cite{PR2-tcs6} for $A\in\{P,D\}$ are
evidently mistaken.

\begin{center}%
\begin{tabular}
[c]{|l|l|l|l|l|l|l|l|l|l|}\hline
$m\setminus x$ & 2 & 3 & 4 & 5 &  & 2 & 3 & 4 & 5\\\hline
100 & 0.50500 & 0.56690 & 0.56429 & 0.56427 &  & 1.33744 & 1.46573 & 1.49445 &
1.51342\\\hline
200 & 0.50250 & 0.56564 & 0.56287 & 0.56286 &  & 1.33203 & 1.46216 & 1.49116 &
1.51038\\\hline
300 & 0.50166 & 0.56522 & 0.56240 & 0.56239 &  & 1.33023 & 1.46097 & 1.49007 &
1.50937\\\hline
400 & 0.50125 & 0.56501 & 0.56216 & 0.56216 &  & 1.32933 & 1.46037 & 1.48952 &
1.50887\\\hline
500 & 0.50100 & 0.56488 & 0.56202 & 0.56202 &  & 1.32879 & 1.46002 & 1.48920 &
1.50856\\\hline
600 & 0.50083 & 0.56480 & 0.56193 & 0.56192 &  & 1.32843 & 1.45978 & 1.48898 &
1.50836\\\hline
700 & 0.50072 & 0.56474 & 0.56186 & 0.56186 &  & 1.32817 & 1.45961 & 1.48882 &
1.50822\\\hline
800 & 0.50063 & 0.56470 & 0.56181 & 0.56181 &  & 1.32798 & 1.45948 & 1.48871 &
1.50811\\\hline
$\vdots$ & $\vdots$ & $\vdots$ & $\vdots$ & $\vdots$ &  & $\vdots$ & $\vdots$
& $\vdots$ & $\vdots$\\\hline
$\infty$ & 0.5 & 0.56438 & 0.56146 & 0.56145 &  & 1.32663 & 1.45860 &
1.48789 & 1.50735\\\hline
\end{tabular}

Table 7A:\ Ratio $m^{a}b_{\left\lfloor x\,m\right\rfloor ,m}/b_{\left\lfloor
x\,m\right\rfloor }$ for Permute ($a=1$) and Graph ($a=1/2$)
\end{center}

\smallskip\ 

\begin{center}%
\begin{tabular}
[c]{|l|l|l|l|l|l|l|l|l|l|}\hline
$m\setminus x$ & 2 & 3 & 4 & 5 &  & 2 & 3 & 4 & 5\\\hline
100 & 0.87413 & 0.95520 & 0.97361 & 0.98570 &  & 1.37273 & 1.54100 & 1.53390 &
1.53385\\\hline
200 & 0.87676 & 0.96022 & 0.97895 & 0.99132 &  & 1.36594 & 1.53756 & 1.53005 &
1.53001\\\hline
300 & 0.87816 & 0.96259 & 0.98148 & 0.99397 &  & 1.36367 & 1.53642 & 1.52876 &
1.52874\\\hline
400 & 0.87906 & 0.96406 & 0.98303 & 0.99559 &  & 1.36254 & 1.53585 & 1.52812 &
1.52810\\\hline
500 & 0.87971 & 0.96508 & 0.98411 & 0.99672 &  & 1.36186 & 1.53551 & 1.52774 &
1.52772\\\hline
600 & 0.88021 & 0.96584 & 0.98492 & 0.99756 &  & 1.36141 & 1.53528 & 1.52748 &
1.52746\\\hline
700 & 0.88060 & 0.96644 & 0.98555 & 0.99823 &  & 1.36109 & 1.53512 & 1.52730 &
1.52728\\\hline
800 & 0.88092 & 0.96693 & 0.98607 & 0.99876 &  & 1.36084 & 1.53500 & 1.52716 &
1.52715\\\hline
$\vdots$ & $\vdots$ & $\vdots$ & $\vdots$ & $\vdots$ &  & $\vdots$ & $\vdots$
& $\vdots$ & $\vdots$\\\hline
$\infty$ & 0.88623 & 0.97438 & 0.99395 & 1.00695 &  & 1.35914 & 1.53415 &
1.52620 & 1.52619\\\hline
\end{tabular}

Table 7B:\ Ratio $m^{a}b_{\left\lfloor x\,m\right\rfloor ,m}/b_{\left\lfloor
x\,m\right\rfloor }$ for Map ($a=1/2$) and Derange ($a=1$)
\end{center}

Why have we devoted effort to evaluating Buchstab's omega? \ Answer: an array
of formulas, parallel to those involving $\rho_{a}(x)$, corresponding to
$h^{\text{th}}$ moments of the smallest component size, were proposed in
\cite{PR2-tcs6}:%
\[%
\begin{array}
[c]{ccccc}%
{\displaystyle\int\limits_{2}^{\infty}}
\dfrac{\omega_{A}(x)}{x^{h+1/2}}dx, &  & h=1,2,\ldots, &  & \text{for }%
A\in\{G,M\}
\end{array}
\]
\ 

\noindent and would be exceedingly attractive. \ Unfortunately the potential
for fulfillment is not good. \ No high-precision numerical estimates of these
integrals are currently known; thus we are not certain that any of the various
$\omega_{A}(x)$ are necessarily allied with $S_{k,n}$ asymptotics in Sections
1--4. \ For now, the formulas remain frustratingly non-explicit and unverified.

\section{Generalized Dickman Rho (II)}

Define $b_{n,m}$ to be the number of $n$-objects whose largest component has
size $\leq m$ (as in Section 5). \ When restricting attention to permutations,
we observe that
\[
\lim_{m\rightarrow\infty}\frac{\left\lfloor x\,m\right\rfloor
\,c_{\left\lfloor x\,m\right\rfloor }}{b_{\left\lfloor x\,m\right\rfloor ,m}%
}=\frac{1}{\rho(x)}%
\]
for any $x>1$, where $\rho(x)$ is the standard Dickman function. \ When
restricting attention to derangements, a factor of $e$ needs to be included
(just as in Sections 4 and 7). \ Replace now the initial factor $\left\lfloor
x\,m\right\rfloor $ in the numerator by $\left\lfloor x\,m\right\rfloor
^{1/2}$. \ Again $2$-regular graphs and mappings possess non-equivalent limits
different from $1/\rho(x)$. \ Just as we found the defining limit for
Buchstab's omega in terms of $\Omega_{a}(x)$ earlier, here we discover the
limit in terms of $\rho_{a}(x)$: \
\[
\lim_{m\rightarrow\infty}\frac{\left\lfloor x\,m\right\rfloor ^{a}%
\,c_{\left\lfloor x\,m\right\rfloor }}{b_{\left\lfloor x\,m\right\rfloor ,m}%
}=\frac{1}{\rho_{a}(x)}\cdot\left\{
\begin{array}
[c]{lll}%
1 &  & \text{if }A=P\text{ and }a=1,\\
e^{3/4}\sqrt{\pi}/2 &  & \text{if }A=G\text{ and }a=1/2,\\
\sqrt{\pi/2} &  & \text{if }A=M\text{ and }a=1/2,\\
e &  & \text{if }A=D\text{ and }a=1.
\end{array}
\right.
\]

We know that $\rho_{a}$ is important (Section 5) and believe that $\omega_{A}$
deserves further study (Section 7). \ It is hoped that someone else might
succeed in carrying on research where we have stopped. \ 

\begin{center}%
\begin{tabular}
[c]{|l|l|l|l|l|l|l|l|l|l|}\hline
$m\setminus x$ & 2 & 3 & 4 & 5 &  & 2 & 3 & 4 & 5\\\hline
100 & 3.23262 & 20.1473 & 196.264 & 2668.39 &  & 15.9879 & 227.490 & 5106.02 &
161434.\\\hline
200 & 3.24569 & 20.3580 & 199.898 & 2742.41 &  & 15.9005 & 227.928 & 5161.08 &
164830.\\\hline
300 & 3.25007 & 20.4291 & 201.130 & 2767.66 &  & 15.8719 & 228.098 & 5180.23 &
166005.\\\hline
400 & 3.25227 & 20.4648 & 201.750 & 2780.41 &  & 15.8577 & 228.188 & 5189.95 &
166600.\\\hline
500 & 3.25359 & 20.4863 & 202.124 & 2788.09 &  & 15.8492 & 228.244 & 5195.83 &
166961.\\\hline
600 & 3.25447 & 20.5006 & 202.373 & 2793.22 &  & 15.8436 & 228.281 & 5199.78 &
167202.\\\hline
700 & 3.25510 & 20.5109 & 202.552 & 2796.90 &  & 15.8396 & 228.308 & 5202.60 &
167375.\\\hline
800 & 3.25558 & 20.5186 & 202.686 & 2799.66 &  & 15.8365 & 228.329 & 5204.73 &
167504.\\\hline
$\vdots$ & $\vdots$ & $\vdots$ & $\vdots$ & $\vdots$ &  & $\vdots$ & $\vdots$
& $\vdots$ & \\\hline
$\infty$ & 3.25889 & 20.5726 & 203.628 & 2819.09 &  & 15.8156 & 228.476 &
5219.73 & 168421.\\\hline
\end{tabular}

Table 8A:\ Ratio $\left\lfloor x\,m\right\rfloor ^{a}c_{\left\lfloor
x\,m\right\rfloor }/b_{\left\lfloor x\,m\right\rfloor ,m}$ for Permute ($a=1$)
and Graph ($a=1/2$)
\end{center}

\smallskip

\begin{center}%
\begin{tabular}
[c]{|l|l|l|l|l|l|l|l|l|l|}\hline
$m\setminus x$ & 2 & 3 & 4 & 5 &  & 2 & 3 & 4 & 5\\\hline
100 & 11.0531 & 165.524 & 3883.12 & 128134. &  & 8.93117 & 55.9354 & 546.943 &
7458.53\\\hline
200 & 10.9698 & 163.012 & 3811.16 & 125549. &  & 8.89477 & 55.9254 & 550.162 &
7558.92\\\hline
300 & 10.9164 & 161.547 & 3766.89 & 123839. &  & 8.88269 & 55.9236 & 551.265 &
7593.21\\\hline
400 & 10.8799 & 160.574 & 3737.04 & 122663. &  & 8.87665 & 55.9229 & 551.822 &
7610.52\\\hline
500 & 10.8531 & 159.869 & 3715.21 & 121794. &  & 8.87304 & 55.9226 & 552.159 &
7620.95\\\hline
600 & 10.8322 & 159.327 & 3698.37 & 121120. &  & 8.87063 & 55.9224 & 552.384 &
7627.93\\\hline
700 & 10.8155 & 158.894 & 3684.85 & 120577. &  & 8.86890 & 55.9223 & 552.545 &
7632.92\\\hline
800 & 10.8016 & 158.537 & 3673.69 & 120127. &  & 8.86761 & 55.9223 & 552.665 &
7636.67\\\hline
$\vdots$ & $\vdots$ & $\vdots$ & $\vdots$ & $\vdots$ &  & $\vdots$ & $\vdots$
& $\vdots$ & \\\hline
$\infty$ & 10.5652 & 152.628 & 3486.92 & 112510. &  & 8.85859 & 55.9221 &
553.517 & 7663.07\\\hline
\end{tabular}

Table 8B:\ Ratio $\left\lfloor x\,m\right\rfloor ^{a}c_{\left\lfloor
x\,m\right\rfloor }/b_{\left\lfloor x\,m\right\rfloor ,m}$ for Map ($a=1/2$)
and Derange ($a=1$)
\end{center}

\section{Addendum}

At the conclusion of Section 5, we gave expressions for the $h^{\text{th}}$
moments of largest component size, given $a=1$ or $a=1/2$, without
justification. \ Here is a plausibility argument. Assuming (absent any proof)
that the first moment for arbitrary $a>0$ is \cite{Lag-tcs6}%
\[
\lambda_{a}=1-%
{\displaystyle\int\limits_{1}^{\infty}}
\frac{\rho_{a}(x)}{x^{2}}dx,
\]
we reverse integration-by-parts:%
\[%
\begin{array}
[c]{lll}%
du=x^{-2-a}dx, &  & v=x^{a}\rho_{a}(x),\\
u=-\dfrac{1}{1+a}x^{-1-a}, &  & dv=a\left[  x^{-1+a}\rho_{a}(x)-(x-1)^{-1+a}%
\rho_{a}(x-1)\right]  dx
\end{array}
\]
and obtain%
\[
-\dfrac{a}{1+a}\left[
{\displaystyle\int\limits_{1}^{\infty}}
\frac{\rho_{a}(x)}{x^{2}}dx-%
{\displaystyle\int\limits_{1}^{\infty}}
\frac{\rho_{a}(x-1)}{(x-1)^{1-a}x^{1+a}}dx\right]  =\left.  -\dfrac{1}%
{1+a}\frac{\rho_{a}(x)}{x}\right\vert _{1}^{\infty}-%
{\displaystyle\int\limits_{1}^{\infty}}
\frac{\rho_{a}(x)}{x^{2}}dx
\]
i.e.,%
\[
\dfrac{a}{1+a}%
{\displaystyle\int\limits_{0}^{\infty}}
\frac{\rho_{a}(x)}{x^{1-a}(x+1)^{1+a}}dx=\dfrac{1}{1+a}-(1-\lambda_{a}%
)+\dfrac{a}{1+a}(1-\lambda_{a})=\dfrac{\lambda_{a}}{1+a}.
\]
From here, we infer (again, absent any proof) that the $h^{\text{th}}$ moment
is
\[%
\begin{array}
[c]{ccc}%
a%
{\displaystyle\int\limits_{0}^{\infty}}
\dfrac{\rho_{a}(x)}{x^{1-a}(x+1)^{h+a}}dx, &  & h=1,2,\ldots\text{.}%
\end{array}
\]

An explanation for $_{S}\mu_{n}$ and $_{S}\sigma_{n}^{2}$ asymptotics in
Sections 2 and 3 was not given. \ Reason: it is more complicated than the
proof in Section 4. \ At some later point, we hope to study combinatorial
objects called \textit{cyclations} \cite{Pip-tcs6} for which moments are known
to be precisely $_{S}G_{1/2}$. \ Since $\sqrt{n}c_{n}/b_{n}\rightarrow
\sqrt{\pi}/2$ as $n\rightarrow\infty$ for these, and because corresponding
limits for graphs and mappings appear in Section 7, the factors $e^{3/4}$ and
$\sqrt{2}$ emerge.

\section{Acknowledgements}

I am grateful to Nicholas Pippenger \cite{Pip-tcs6} for a helpful discussion.
\ The creators of Mathematica, as well as administrators of the MIT Engaging
Cluster, earn my gratitude every day. \ A sequel to this paper is now
available \cite{F1-tcs6}; another will be released soon \cite{F2-tcs6}.

\end{document}